\documentclass[10pt]{article}
\usepackage{geometry}                % See geometry.pdf to learn the layout options. There are lots.
\usepackage{graphicx}
\usepackage{amssymb}
\usepackage{amsmath}
\usepackage{epstopdf}
\usepackage{stmaryrd}
\usepackage{bm}
\usepackage{accents}
\usepackage{adjustbox}
\usepackage{stmaryrd}

\usepackage{color}
\definecolor{redcol}{rgb}{1.,0.,0.0} 
\definecolor{lnkcol}{rgb}{0.,0.,0.0} %black

\usepackage[pdftex,pdfpagelabels=true]{hyperref}
\hypersetup{  colorlinks=true,
              citecolor=lnkcol,%
              linkcolor=lnkcol,%
              urlcolor=lnkcol,
            pdfborder=false,
            bookmarks=true,
            bookmarksnumbered=true,
            pdftitle = {working title: BR1 stability},
pdfsubject = {},
pdfauthor = {},
pdfkeywords = {}
 }
\usepackage[all]{hypcap}  % hyperlink soll oberhalb des Bilds anhalten
\usepackage[dvips]{thumbpdf}

\usepackage{cleveref} % === > !!! MUST BE USED AFTER hyperref!!! <====
\usepackage{stackengine}
\stackMath

\usepackage{accsupp} % for ensuring the right Unicode codepoint upon pasting
\newcommand*{\llbrace}{%
  \BeginAccSupp{method=hex,unicode,ActualText=2983}%
    \textnormal{\usefont{OMS}{lmr}{m}{n}\char102}%
    \hspace*{-3pt}
    \textnormal{\usefont{OMS}{lmr}{m}{n}\char102}%
  \EndAccSupp{}%
}
\newcommand*{\rrbrace}{%
  \BeginAccSupp{method=hex,unicode,ActualText=2984}%
    \textnormal{\usefont{OMS}{lmr}{m}{n}\char103}%
        \hspace*{-3pt}
    \textnormal{\usefont{OMS}{lmr}{m}{n}\char103}%
  \EndAccSupp{}%
}
\makeatletter
\def\munderbar#1{\underline{\sbox\tw@{$#1$}\dp\tw@\z@\box\tw@}}
\makeatother

\newcommand\iprod[1]{\left\langle #1\right\rangle} 				% inner product
\newcommand\inorm[1]{\left |\left| #1\right|\right|}		% norm
\newcommand\iprodN[1]{\left\langle #1\right\rangle_{N}}			% discrete inner product subscripted with N
\newcommand\inormN[1]{\left |\left| #1\right|\right|_{N}}	% discrete norm subscripted with N
\newcommand\spacevec[1]{\accentset{\,\rightarrow}{#1}}					% spatial vector, i.e. x \hat x + y \hat y
\newcommand\contravec[1]{\tilde{ #1}}					% contravariant vector
\newcommand\statevec[1]{\mathbf #1}					% state vector, e.g. [rho rhou rhoe]^{T}

\newcommand\dS{\,\operatorname{dS} }

\newcommand\PBT{\,\operatorname{PBT} }

\newcommand\mmatrix[1]{\underbar{#1}}				% Matrix as taught in linear algebra, i.e. math matrix.
\newcommand{\jump}[1]{\left\llbracket#1\right\rrbracket}   % jump at element boundary
\newcommand{\average}[1]{\left\{\!\!\left\{#1\right\}\!\!\right\}}   % average at element boundary
\newcommand{\partialderivative}[2]{\frac{\partial #1}{\partial #2}}  % to make first partial derivatives
\newcommand{\avg}[1]{\ensuremath{\llbrace#1\rrbrace}}
\newcommand{\inPN}{\in\mathbb{P}^{N}}
\newcommand{\IN}[1]{\mathbb{I}^{N}\left( #1\right)}

\newcommand\interiorfaces{{\mathrm{interior}\atop\mathrm{faces}}}
\newcommand\boundaryfaces{{\mathrm{Boundary}\atop\mathrm{faces}}}

\title{Stability of Overintegration Methods for Nodal Discontinuous Galerkin Spectral Element Methods}
\author{David A. Kopriva}
\date{Received: date / Accepted: date}

\begin{document}
\maketitle
\begin{abstract}
We perform stability analyses for discontinuous Galerkin spectral element approximations of linear variable coefficient hyperbolic systems in three dimensional domains with curved elements. Although high order, the precision of the quadratures used are typically too low with respect to polynomial order associated with their arguments, which introduces aliasing errors that can destabilize an approximation, especially when the solution is underresolved. We show that using a larger number of points in the volume quadrature, often called ``overintegration'', can eliminate the aliasing term associated with the volume, but introduces new aliasing errors at the surfaces that can destabilize the solution. Increased quadrature precision on both the volume and surface terms, on the other hand, leads to a stable approximation. The results support the findings of Mengaldo \textit{et al.} [\textit{Dealiasing techniques for high-order spectral element methods on
  regular and irregular grids. Journal of Computational Physics, 299:56 -- 81, 2015}] who found that fully consistent integration was more robust for the solution of compressible flows than the volume only version.
\end{abstract}
\section{Introduction}
Discontinuous Galerkin (DG) spectral element methods (DGSEM) for hyperbolic systems of equations are characterized by an approximation of a weak form of the equations where the solutions and fluxes are approximated by polynomials of degree $N$ (typically large, $N\ge 4$) with nodal values represented at Gauss (or Gauss-Lobatto) points, and inner products are approximated by the associated Gauss quadrature. Unfortunately, the volume and surface flux terms usually are underintegrated with respect to the polynomial order of the advective fluxes. This happens in nonlinear problems, but 
also for variable coefficient linear problems, or even constant coefficient problems posed on curved element meshes. 

Although the method usually works in practice, aliasing errors associated with the underintegration of the fluxes can lead to instability, especially if the solution is severely underresolved \cite{Beck:2013sf},\cite{Gassner:2013qf}.  In \cite{Kirby2003}, the authors proposed a method to stabilize nodal DGSEMs for hyperbolic systems, which they called ``super-collocation'' and which, though perhaps better called ``consistent integration'', is commonly called ``overintegration''. In  \cite{Kirby2003}, Gauss quadratures of sufficiently high order were used to approximate the volume integrals so as to eliminate the polynomial aliasing associated with the quadrature. Later, in \cite{Mengaldo201556}, numerical experiments suggested that overintegrating the surface integrals increases the robustness of approximations of the Euler gas dynamics equations. Overintegration is used in large codes such as Nektar++\cite{ca.mo.co.ea:15}
%\cite{Cantwella:2015hs}
and Flexi\cite{Hindenlang201286} for stabilization. However theoretical justification of the procedures has not formally been established to date.

In this paper we study the linear (energy) stability for hyperbolic systems with variable coefficients and for curved elements in three space dimensions of (i) the original DG spectral element approximation, (ii) overintegration of the volume terms only, and (iii) overintegration of the surface and volume terms together. Our results support the conclusion of \cite{Mengaldo201556} that overintegration of both the surface and volume terms is more robust than overintegration of the volume terms alone.

\section{The Hyperbolic Boundary-Value Problem}
In this paper we study the stability of discontinuous Galerkin spectral element methods for boundary-value problems for variable coefficient linear hyperbolic systems of the form 
\begin{equation}
{{\mathbf{u}}_t} + \nabla  \cdot \spacevec {\mathbf{f}}  = 0
\label{eq:ConservativeSystem}
\end{equation}
defined on a domain $\Omega$. Here, $\mathbf{u}\left(\spacevec x,t\right) = [u_{1}\left(\spacevec x,t\right)\;u_{2}\left(\spacevec x,t\right)\;\ldots \; u_{p}\left(\spacevec x,t\right)]^{T}$ is the state vector and 
\begin{equation}
\spacevec {\mathbf{f}} \left(\mathbf u\right) = \sum\limits_{m = 1}^3 {{\mmatrix{A}_{m}}\left( {\spacevec x} \right){\mathbf{u}}{{\hat x}_m}}  = \spacevec {\mmatrix{A}}\left(\spacevec x\right) {\mathbf{u}}
\end{equation}
is the linear flux space vector, where $\mmatrix A_{m}$ is the coefficient matrix of the derivative in the $m^{th}$ space dimension and $\hat x_{m}$ is a unit vector in that direction. We assume that the system has already been symmetrized and is hyperbolic, that is
\begin{equation}
\mmatrix{A}_{m} = \left(\mmatrix{A}_{m}\right)^{T}\quad\text{and}\quad \sum\limits_{m = 1}^3 {{\alpha _m}{\mmatrix{A}_{m}}}  = \mmatrix{R}\Lambda {\mmatrix{R}^{ - 1}}
\label{eq:mDHyperbolicConditions_DAK}
\end{equation}
for any $\sum\limits_{m = 1}^3 {\alpha _m^2}  \ne 0$ and some real diagonal matrix $\Lambda$. We  also assume that the matrices  are continuous in their argument and that the $\mmatrix{A}_{m}$ have bounded derivatives in the sense that 
\begin{equation}{\left\| {\nabla  \cdot \spacevec{ \mmatrix{A}} } \right\|_2}<\infty,\end{equation}
where $\left\|\cdot\right\|_2$ is the matrix 2-norm. 
Adding the initial condition $\statevec u (\spacevec x, 0) = \statevec u_{0}(\spacevec x)$ and external state boundary values $\statevec g$ properly imposed  on $\partial\Omega$ completes the specification of the boundary-value problem.

With characteristic boundary conditions specified along the $\partial \Omega$, the $\mathbb L^{2}$ norm of the solution  $\left\| {\mathbf{u}} \right\| = \sqrt {\left\langle {{\mathbf{u}},{\mathbf{u}}} \right\rangle}$ defined through the inner product $
\left\langle {{\mathbf{u}},{\mathbf{v}}} \right\rangle = \int_\Omega  {{{\mathbf{u}}^T}{\mathbf{v}}dxdydz}
$
is bounded in terms of the initial and boundary data as \cite{Kopriva:2017yg}
\begin{equation}
\frac{d}{{dt}}{\left\| \statevec u \right\|^2} =  - \int_{\partial \Omega } {\left\{ {{\statevec u^T}{\mmatrix A^ + }\statevec u - {\statevec g^T}\left| {{\mmatrix A^ - }} \right|\statevec g} \right\}dS}  + 2\gamma {\left\| u \right\|^2},
\label{eq:TrueenergyTimeDerivative}
\end{equation}
where $\gamma = \frac{1}{2}\mathop {\max }\limits_\Omega  \left|\left| \nabla\cdot\spacevec{\mmatrix {A}} \right|\right|_{2}$, 
\begin{equation}
\spacevec{\mmatrix{A}} \cdot \hat n = \sum\limits_{m = 1}^3 {{\mmatrix A_{m}}{{\hat n}_m}}  = \mmatrix R\Lambda {\mmatrix R^{ - 1}} = \mmatrix R{\Lambda ^ + }{\mmatrix R^{ - 1}} + \mmatrix R{\Lambda ^ - }{\mmatrix R^{ - 1}} \equiv {\mmatrix A^ + } + {\mmatrix A^ - },
\end{equation}
$\Lambda^{\pm}=\left(\Lambda \pm \left|\Lambda\right|\right)/2$ are diagonal matrices, and $\hat n$ is the outward normal to the surface of the domain.
Integrating \eqref{eq:TrueenergyTimeDerivative} in time gives
\begin{equation}
\begin{split}
\inorm{\statevec{u}(T)}^{2} +\int_{0}^{T} {\int_{\partial \Omega } {{{\statevec{u}}^T}{\mmatrix{A}^{+}}{\statevec u}dS}dt}
&\le e^{2\gamma T}\inorm{\statevec{u}(0)}^{2}+\int_{0}^{T} {\int_{\partial \Omega } {e^{2\gamma(T-t)}{{\statevec{g}}^T}\left|\mmatrix{A}^{-}\right|{\statevec g}dS}dt} 
\\&\le e^{2\gamma T}\left\{\inorm{\statevec{u}_{0}}^{2}+\int_{0}^{T} {\int_{\partial \Omega } {{{\statevec{g}}^T}\left|\mmatrix{A}^{-}\right|{\statevec g}dS}dt}\right\}.
\end{split}
\label{eq:dDimensionWellPosednessLinear_DAK}
\end{equation}

Equations (\ref{eq:TrueenergyTimeDerivative}) and (\ref{eq:dDimensionWellPosednessLinear_DAK}) show that energy is dissipated at the physical boundary along outgoing characteristics ($\mmatrix A^{+}$ term), added through the incoming characteristics ($\mmatrix A^{-}$ term), and can grow as a result of the variable coefficients in the problem ($e^{2\gamma t}$ factor).
If the external state $\statevec g=0$ and $\nabla\cdot\spacevec{\mmatrix A} = 0$ then $\gamma = 0$ and (\ref{eq:dDimensionWellPosednessLinear_DAK}) reduces to
\begin{equation}
\inorm{\statevec u(T)}\le\inorm{\statevec u_{0}},
\end{equation}
so that the energy does not grow.

\section{DG Spectral Element Methods for Variable Coefficient Problems on Curved Elements}
We consider here the domain $\Omega$ subdivided into $N_{el}$ nonoverlapping hexahedral elements, $e^{k}, k=1,2,\ldots,N_{el}$. We assume that the subdivision is conforming. Each element is mapped from the reference element $E$ by a transformation $\spacevec x = \spacevec X\left(\spacevec \xi\right)$. 
Directions in physical space can be represented in terms of the mapping by the three covariant basis vectors
\begin{equation}
\spacevec{a}_{i}=\partialderivative{\spacevec{X}}{\xi^{i}}\quad i=1,2,3,
\end{equation}
and the (volume weighted) contravariant vectors, formally written as
\begin{equation}
\mathcal{J}\spacevec{a}^{i}=\spacevec{a}_{j}\times\spacevec{a}_{k}, \quad (i,j,k)\;\text{cyclic},
\end{equation}
where
\begin{equation}
\mathcal{J}=\spacevec{a}_{1}\cdot\left(\spacevec{a}_{2}\times \spacevec{a}_{3}\right)
\end{equation}
is the Jacobian of the transformation. 
We will assume that the metric terms are approximated as polynomials of degree $N$ so that
\begin{equation}
\sum\limits_{i = 1}^3 {\frac{{\partial {\mathbb{I}^N}\left( {J{{\spacevec a}^i}} \right)}}{{\partial {\xi ^i}}}}  = 0,
\end{equation}
as described in \cite{kopriva2006},\cite{Kopriva2016274}, where $\mathbb{I}^{N}:L^{2}(E)\rightarrow\mathbb{P}^{N}$ is the polynomial interpolation operator and $\mathbb{P}^{N}$ is the space of polynomials of degree less than or equal to $N$ on $E$.

Under the mapping, the 
divergence of a spatial vector flux can be written compactly in terms of the reference space variables
as
\begin{equation}
\nabla \cdot \spacevec {\statevec f} 
= \frac{1}{\mathcal{J}}\sum\limits_{i = 1}^3 {\frac{\partial }{{\partial {\xi ^i}}}\left( {\mathcal{J}{{\spacevec a}^i} \cdot \spacevec{ \statevec f}} \right)}
= \frac{1}{\mathcal{J}}\sum\limits_{i = 1}^3 {\frac{\partial \tilde {\statevec f}^{i}}{{\partial {\xi ^i}}}}
= \frac{1}{\mathcal{J}}{\nabla _\xi } \cdot \contravec{ {\statevec f}}.
\label{eq:CompSpaceDivergence_DAK}
\end{equation}
The vector $\contravec{{\statevec f}}$ is the volume weighted contravariant flux whose components are $\contravec{{\statevec f}}^{i}={\mathcal{J}{{\spacevec a}^i} \cdot \spacevec{ {\statevec f}}}$.
The conservation law can therefore be represented on the reference domain by another conservation law
\begin{equation}
{\mathcal{J}\statevec u_t} +\nabla_{\xi}\cdot\left(\contravec{\mmatrix{A}}\statevec u\right)  = 0,
\label{eq:mappedEquation}
\end{equation}
where we define the (volume weighted) contravariant coefficient matrices
\begin{equation}
\contravec{\mmatrix A}^i = \mathcal{J}{{\spacevec a}^i} \cdot \spacevec{\mmatrix{A}}
\end{equation}
and
\begin{equation}\contravec{\mmatrix{A}} = \sum\limits_{i = 1}^3 {\contravec{\mmatrix A}^i{{\hat \xi }^i}}. \end{equation}

The discontinuous Galerkin approximation is created from a weak formulation formed by multiplying (\ref{eq:mappedEquation}) by a test function $\boldsymbol \phi$, forming the inner product over the reference element, $E$, and integrating by parts.
\begin{equation}
\iprod{\mathcal J\statevec u,\boldsymbol\phi}+\int_{\partial E}{\boldsymbol\phi^{T}\contravec{\statevec f}\cdot\hat ndS}-\iprod{\contravec{\statevec f},\nabla\boldsymbol\phi} = 0.
\label{eq:WeakContinuousForm}
\end{equation}

The solution and fluxes are then approximated by polynomials global within each element.
The solution $\statevec u$ is approximated by a polynomial of degree $N$, $\statevec u \approx \statevec U\in\mathbb{P}^{N}$ as is the Jacobian, $\mathcal J\approx J\inPN$. They are polynomial interpolants with nodes at the nodes of the Gauss quadrature used later to approximate the inner products. 

With variable coefficients and variable metric terms, the fluxes are not necessarily polynomials of degree $N$ or less even if the solution is a polynomial of degree $N$. They can be approximated as polynomials of degree $N$ in several ways. One can, for instance, approximate the contravariant flux as a polynomial
of degree $2N$
\begin{equation}
 \contravec{\statevec F} = \contravec{\mathcal A}\statevec U\in\mathbb{P}^{2N}.
 \label{eq:fInP2N}
\end{equation}
by approximating the coefficient matrices as a polynomial of degree $N$
\begin{equation}
\contravec{\mathcal A}= \mathbb{I}^{N}\left( \mathbb{I}^{N}\left(J\spacevec a^{i}\right)\cdot \spacevec{\mmatrix A} \right)=\mathbb{I}^{N}\left( \mathbb{I}^{N}\left(J\spacevec a^{i}\right)\cdot \mathbb I^{N}\left(\spacevec{\mmatrix A}\right) \right)\inPN.
\label{eq:CoeffMatrixPoly}
\end{equation}
Alternatively, one could approximate the contravariant coefficient matrices
by
\begin{equation}
\contravec{\mathcal A}=  \mathbb{I}^{N}\left(J\spacevec a^{i}\right)\cdot \mathbb I^{N}\left(\spacevec{\mmatrix A}\right) \in \mathbb{P}^{2N},
\label{eq:CoeffMatrixPoly2N}
\end{equation}
leading to a contravariant
flux
\begin{equation}
\contravec{\statevec F} = \contravec{\mathcal A}\statevec U\in\mathbb{P}^{3N}.
 \label{eq:fInP3N}
\end{equation}
Finally, the usual (and underintegrated) approximation interpolates the entire flux, 
\begin{equation}
 \contravec{\statevec F} = \mathbb{I}^{N}\left( \mathbb{I}^{N}\left(J\spacevec a^{i}\right)\cdot \spacevec{\mmatrix A} \statevec U\right)\in\mathbb{P}^{N}.
  \label{eq:fInPN}
\end{equation}
If none of these is done, and $\spacevec{\mmatrix A}$ is not a polynomial function of the reference space variables, then the flux $ \contravec{\statevec F} $ is not a polynomial function approximation. 

In the following, we will assume that
\begin{equation}
{\mathbf{U}} \in {\mathbb{P}^N},\quad \tilde {\mathbf{F}}  \in {\mathbb{P}^{pN}},
\end{equation}
for some $p$.
As an important aside, the fluxes of the Euler equations of gas dynamics are not polynomials of the state vector of the conservative variables, and hence there is no exact representation of the fluxes as a polynomial when the state vector is a polynomial \cite{Gassner:2016ye}.

The normal surface flux in \eqref{eq:WeakContinuousForm} is replaced by an interpolant, $\tilde{\statevec F}^{*}$, of a numerical flux function, $\tilde{\statevec f}^{*}\left(\statevec U^{L},\statevec U^{R}\right)$, that is a function of the left and right states at the faces and is continuous across the face. The exact upwind numerical flux for the linear flux function is
\begin{equation}
\tilde{\statevec f}^{*}\left(\statevec U^{L},\statevec U^{R}\right)=\frac{\tilde{\mmatrix A}\cdot\hat n\statevec U^{L} + \tilde{\mmatrix A}\cdot\hat n\statevec U^{R}}{2}-\frac{1}{2} \left| \tilde{\mmatrix A}\cdot\hat n\right|\left(\statevec U^{R} - \statevec U^{L}\right) \equiv \tilde{\mmatrix A}\cdot\hat n\avg{\statevec U} - \frac{1}{2}\left| \tilde{\mmatrix A}\cdot\hat n\right|\jump{\statevec U},
\label{eq:ExactNumericalFlux}
\end{equation}
where $\hat n$ is the reference space normal and the usual jump and average operators are $\jump{\cdot}$ and $\average{\cdot}$, respectively.

The final equality in \eqref{eq:ExactNumericalFlux} will be true for the discrete approximation if the interpolant of the coefficient matrices, $\tilde{\mmatrix A}$, is continuous at the interfaces.  We can ensure continuity if the interpolation points include the boundary points, as would be the case if the Gauss-Lobatto nodes are used as the interpolation points, but not if the Gauss points are used. 
For this reason, and the fact that they were used in \cite{Mengaldo201556}, we consider only approximations that use the Gauss-Lobatto points in this paper.
For simplicity, we define  $\tilde {\mmatrix A}^{\pm} = \frac{1}{2}\left(\tilde {\mmatrix A}\cdot\hat n \pm \left|\tilde {\mmatrix A}\cdot\hat n\right| \right)$ so that $\tilde{\statevec f}^{*}\left(\statevec U^{L},\statevec U^{R}\right) = \tilde {\mmatrix A}^{+}\statevec U^{L}+\tilde {\mmatrix A}^{-}\statevec U^{R}.$ As with the interior fluxes, \eqref{eq:fInP2N}, \eqref{eq:fInP3N}, \eqref {eq:fInPN},
the numerical surface flux would  be a polynomial of degree $3N$, $2N$ or $N$ depending on how it is approximated. 

Inserting the polynomial approximations and the numerical surface flux gives us the exactly integrated discontinuous Galerkin approximation
\begin{equation}\left\langle {{{\mathbf{U}}_t},\phi } \right\rangle  + \int_{\partial E} {{{\mathbf{F}}^{*,T}}\phi dS}  - \left\langle {{\mathbf{\spacevec F}},\nabla \phi } \right\rangle  = 0,
\label{eq:approxWExactIntegration}
\end{equation}
where ${{\mathbf{F}}^*}, {\mathbf{\spacevec F}} \in {\mathbb{P}^{pN}}$.

Because the arguments are high order polynomials, the inner products are not integrated exactly in practice. Instead, they are approximated by quadrature. We write the Gauss-Lobatto quadrature over $E$ as
 \begin{equation}
 \int_{E,N} {gd \xi d\eta d\zeta}  \equiv \sum\limits_{i,j,k = 0}^N {{g_{ijk}}{\omega_{i}}{\omega_{j}}{\omega_{k}}},
  \end{equation}
where $\omega_{i}$, $\omega_{j}$ and $\omega_{k}$ are the quadrature weights for the $i,j,$ and $k$ coordinate directions, and $g_{ijk} = g\left(\xi_{i},\eta_{j},\zeta_{k}\right)$ are the values of $g$ evaluated at the quadrature points.
Two-dimensional surface integral approximations of a vector function $\spacevec g$ are
 \begin{equation}
 \begin{split}
 \int_{\partial E,N} {\spacevec g \cdot \hat ndS}  &\equiv \sum\limits_{i,j = 0}^N {\left. {{\omega_{ij}}{g^{(1)}}\left( {\xi,{\eta_i},{\zeta_j}} \right)} \right|_{\xi =  - 1}^1}  + \sum\limits_{i,j = 0}^N {\left. {{\omega_{ij}}{g^{(2)}}\left( {{\xi_i},\eta,{\zeta_j}} \right)} \right|_{\eta =  - 1}^1}  + \sum\limits_{i,j = 0}^N {\left. {{\omega_{ij}}{g^{(3)}}\left( {{\xi_i},{\eta_j},\zeta} \right)} \right|_{\zeta =  - 1}^1} 
\\&\equiv\int_N {\left. {{g^{(1)}}d\eta d\zeta } \right|} _{\xi  =  - 1}^1 + \int_N {\left. {{g^{(2)}}d\xi d\zeta } \right|} _{\eta  =  - 1}^1 + \int_N {\left. {{g^{(3)}}d\xi d\eta } \right|} _{\zeta  =  - 1}^1,
 \end{split}
 \end{equation}
 where $\omega_{ij} \equiv \omega_{i}\omega_{j}$, etc.
 Finally, we define the discrete inner product of two functions $\statevec f$ and $\statevec g$ and the discrete norm of $\statevec f$ from the quadrature
 \begin{equation}{\iprod{\statevec f,\statevec g}_{E,N}} = \int_{E,N} {\statevec f^{T} \statevec gd \xi d \eta d \zeta }  \equiv \sum\limits_{i,j,k = 0}^N {{\statevec f^{T}_{ijk}}{\statevec g_{ijk}}{\omega_{ijk}}} ,\quad {\left\| \statevec f \right\|_{E,N}} = \sqrt {{\iprod{\statevec f,\statevec f}_{E,N}}}. \end{equation}
 We often use the fact that the definition of the discrete inner product implies that
 \begin{equation}
 \iprod{\statevec f,\statevec V}_{E,N} =  \iprod{\mathbb{I}^{N}(\statevec f),\statevec V}_{E,N}\quad\forall \statevec V\inPN.
 \label{eq:DL2Projection}
 \end{equation}
 
 Stability analysis hangs on the fact that the Gauss-Lobatto quadrature rule satisfies a summation by parts property \cite{gassner2010} and from that, the discrete Gauss law \cite{Kopriva:2017yg}
  \begin{equation}
 {\iprod{\nabla  \cdot \spacevec{\statevec F},\statevec V}_{E,N} = \int_{\partial E ,N} {\statevec V^{T}\spacevec {\statevec F} \cdot \hat ndS}  - \iprod{\spacevec {\statevec F},\nabla \statevec V}_{E,N}},\;\forall\; \statevec V\inPN.
 \label{eq:DiscreteGreens_DAK}
 \end{equation}

With the definitions for the quadrature, we write a general DG spectral element approximation for (\ref{eq:approxWExactIntegration}), where the inner products are approximated by possibly different polynomial orders $N,M,L$ by restricting $\boldsymbol\phi\inPN$ and writing
\begin{equation}
\iprod{ J\statevec U,\boldsymbol\phi}_{N}+\int_{\partial E, L}{\boldsymbol\phi^{T}\contravec{\statevec F}^{*}dS}-\iprod{\contravec{\statevec F},\nabla\boldsymbol\phi}_{M} = 0,
\label{eq:WeakContinuousForm2}
\end{equation}
where we have left off the subscript $E$ to be understood in context. The form (\ref{eq:WeakContinuousForm2}) is generally known as the ``weak'' form of the approximation.

 \section{Stability Analysis}
We now study the stability of three approximations. The first is the standard approximation with $N=L=M$. We show that instability can be caused by large aliasing errors in severely underresolved situations. The second is fully integrated, with $L = M > N$ chosen so that quadrature is exact and the approximation is stable. Finally, we examine the approximation where only the volume term is fully integrated, $L = N, M>N$ to see that for severely underresolved problems, interpolation error along the element surfaces can lead to instability.
\subsection{Standard Approximation: $L=M=N$}
\label{sec:StandardApprox}
The standard approximation with $L=M=N$ is equivalent to approximating $\contravec{\statevec{F}}=\IN{\tilde{\mathcal A}\statevec U}\inPN$ using \eqref{eq:fInPN} by way of \eqref{eq:DL2Projection}.
Using quadrature of the same order for all of the inner products cannot be shown to be stable unless $\contravec{\statevec F} = \contravec{\mathcal A}\statevec U\inPN$ and  $\contravec{\mathcal A}$ is a constant, which corresponds to constant coefficients and rectangular elements. When we set $\boldsymbol\phi = \statevec U$ in (\ref{eq:WeakContinuousForm2}),
\begin{equation}
\iprod{ J\statevec U,\statevec U}_{N}+\int_{\partial E, N}{\statevec U^{T}\contravec{\statevec F}^{*}dS}-\iprod{\contravec{\statevec F},\nabla\statevec U}_{N} = 0,
\label{eq:WeakContinuousFormNNN1}
\end{equation}
or writing the first term as the time derivative of the $J$-weighted norm,
\begin{equation}
\frac{1}{2}\frac{d}{dt}\inorm{\statevec U}^{2}_{J,N}+\int_{\partial E, N}{\statevec U^{T}\contravec{\statevec F}^{*}dS}-\iprod{\contravec{\statevec F},\nabla\statevec U}_{N} = 0.
\label{eq:WeakContinuousFormNNN}
\end{equation}
We then use the discrete Gauss law \eqref{eq:DiscreteGreens_DAK} on the volume term
\begin{equation}
\iprodN{\contravec{\statevec F},\nabla\statevec U}+\iprodN{\nabla\cdot\contravec{\statevec F},\statevec U}=\int_{\partial E, N}{\statevec U^{T}\contravec{\statevec F}\cdot\hat ndS}
\end{equation}
to re-write \eqref{eq:WeakContinuousFormNNN} as
\begin{equation}
\frac{1}{2}\frac{d}{dt}\inorm{\statevec U}^{2}_{J,N}+\int_{\partial E, N}{\statevec U^{T}\left\{\contravec{\statevec F}^{*} -\contravec{\statevec F}\cdot\hat n\right\} dS}+\iprod{\nabla\cdot\contravec{\statevec F},\statevec U}_{N} = 0.
\label{eq:WeakContinuousFormNNNStrong}
\end{equation}

To account for the fact that the product rule does not in general hold for the interpolant of a product, we add and subtract
\begin{equation}
\frac{1}{2}\iprodN{\contravec{\mathcal A}\cdot\nabla\statevec U + \left(\nabla\cdot\contravec{\mathcal A}\right)\statevec U,\statevec U },
\end{equation}
where $\tilde{\mathcal{A}}=\mathbb{I}^{N}\left( \mathbb{I}^{N}\left(J\spacevec a^{i}\right)\cdot \spacevec{\mmatrix A}\right)\inPN$ to write
\begin{equation}
\begin{split}
\iprodN{\nabla\cdot\contravec{\statevec F}, \statevec U} &= 
\frac{1}{2}\iprodN{\nabla\cdot\contravec{\statevec F}, \statevec U} 
+ \frac{1}{2}\iprodN{\contravec{\mathcal A}\cdot\nabla\statevec U,\statevec U}
+\frac{1}{2}\iprodN{\nabla\cdot\contravec{\mathcal A}\statevec U,\statevec U}
\\&-\frac{1}{2}\iprodN{\contravec{\mathcal A}\cdot\nabla\statevec U + \left(\nabla\cdot\contravec{\mathcal A}\right)\statevec U-\nabla\cdot\contravec{\statevec F},\statevec U}.
\end{split}
\end{equation}
Now the discrete Gauss law \eqref{eq:DiscreteGreens_DAK}, the collocation of the solution and flux, and symmetry of $\contravec{\mathcal A}$ imply that
\begin{equation}
\begin{split}
\iprodN{\contravec{\mathcal A}\cdot\nabla\statevec U,\statevec U} &= \int_{\partial E, N}{\statevec U^{T}\tilde{\statevec F}dS}-\iprodN{\statevec U,\nabla\cdot\IN{\tilde{\mathcal A}\statevec U}}
\\&  = \int_{\partial E, N}{\statevec U^{T}\tilde{\statevec F}dS}-\iprodN{\statevec U,\nabla\cdot\tilde{\statevec F}}.
\end{split}
\label{eq:DGLOnFlux}
\end{equation}
Therefore,
\begin{equation}
\begin{split}
\iprod{\contravec{\statevec F},\nabla\statevec U}_{N}
&=
\frac{1}{2}\int_{\partial E, N}{\statevec U^{T}\contravec{\statevec F}\cdot\hat ndS}
-\frac{1}{2}\iprodN{\left(\nabla\cdot\contravec{\mathcal A}\right)\statevec U,\statevec U}
\\&+\frac{1}{2}\iprodN{\contravec{\mathcal A}\cdot\nabla\statevec U 
+ \left(\nabla\cdot\contravec{\mathcal A}\right)\statevec U
-\nabla\cdot\tilde{\statevec F},\statevec U}.
\end{split}
\end{equation}
Then (\ref{eq:WeakContinuousFormNNN}) becomes
\begin{equation}
\begin{split}
\frac{1}{2}\frac{d}{dt}\inorm{\statevec U}^{2}_{J,N}
=&-\int_{\partial E, N}{\statevec U^{T}\left\{\contravec{\statevec F}^{*}-\frac{1}{2}\contravec{\statevec F}\cdot\hat n\right\}dS}
\\&-\frac{1}{2}\iprodN{\left(\nabla\cdot\contravec{\mathcal A}\right)\statevec U,\statevec U}
\\&+\frac{1}{2}\iprodN{\contravec{\mathcal A}\cdot\nabla\statevec U 
+ \left(\nabla\cdot\contravec{\mathcal A}\right)\statevec U
-\nabla\cdot\IN{\contravec{\statevec F}},\statevec U}.
\end{split}
\label{eq:WeakContinuousFormNNN2}
\end{equation}

The last inner product in (\ref{eq:WeakContinuousFormNNN2}) is the projection of the amount by which the product rule is not satisfied by polynomial functions. Note that if $\contravec{\mathcal A}$ is constant, which happens when the original PDE is constant coefficient and the elements are rectangular, then $\contravec{\statevec F}=\IN{\tilde{\mathcal A}\statevec U} = \tilde{\mathcal A}\statevec U$ and the last term vanishes. Otherwise the last term can be positive or negative. 

The product rule error term can be bounded as
\begin{equation}
\iprodN{\contravec{\mathcal A}\cdot\nabla\statevec U 
+ \left(\nabla\cdot\contravec{\mathcal A}\right)\statevec U
-\nabla\cdot\IN{\contravec{\statevec F}},\statevec U}\le2\varepsilon_{A}\inorm{U}_{J,N}^{2},
\end{equation}
where
\begin{equation}
\varepsilon_{A}=\frac{1}{2}\mathop {\max }\limits_{E,||\statevec U||{_{J,N}} \ne 0} \frac{{{{\left\| {\contravec{\mathcal A} \cdot \nabla \statevec U + \left( {\nabla  \cdot \contravec{\mathcal A}} \right)\statevec U - \nabla  \cdot \IN{\contravec{\statevec F}}} \right\|}^{2}_2}}}{J{||\statevec U||^{2}_{J,N}}}.
\end{equation}
Note also that the sign of $\nabla\cdot\contravec{\mathcal A}$ can be positive or negative and that
\begin{equation} 
- \left\langle {\left( {\nabla  \cdot \contravec{\mathcal A}} \right)\statevec U,\statevec U} \right\rangle_{N}  \leqslant 2\hat \gamma \inorm{\statevec U}_{J,N}
\end{equation}
where
\begin{equation}
\hat \gamma = \mathop {\max }\limits_E \frac{{\left\| {\nabla  \cdot \contravec{\mathcal A}} \right\|}_{2}}{J}.
\end{equation}
Therefore, we can bound the elemental energy in \eqref{eq:WeakContinuousFormNNN2} as
\begin{equation}
\frac{d}{dt}\inorm{\statevec U}^{2}_{J,N}\le -2\int_{\partial E, N}{\statevec U^{T}\left\{\contravec{\statevec F}^{*}
-\frac{1}{2}\contravec{\statevec F}\cdot\hat n\right\}dS}
+2\left(\hat\gamma+\varepsilon_{A}\right)\inorm{\statevec U}^{2}_{J,N}.
\label{eq:LocalDiscreteEnergyDGSEM}
\end{equation}

Equation (\ref{eq:LocalDiscreteEnergyDGSEM}) shows how the energy changes on a single element. Summing over all elements $e^{k}$ gives the total energy change
\begin{equation}
\begin{split}
\frac{d}{dt} \sum\limits_{k=1}^{N_{el}}\inorm{\statevec U^{k}}^{2}_{J,N}\le& 2
 \sum\limits_{Interior\atop Faces} \int_{\partial E,N } {\left\{{\statevec F}^{*,T}\jump{\statevec U}-\frac{1}{2}\jump{\left(\contravec {\statevec F}\cdot\hat n\right)^{T}\statevec U} \right\}\dS}- \PBT\\&+ 2\left(\hat\gamma+\varepsilon_{A}\right)\sum\limits_{n=1}^{N_{el}}\inorm{\statevec U^{k}}^{2}_{J,N},
\end{split}
\label{eq:GlobalStandardEnergyTimeDerivative}
\end{equation}
where $\hat\gamma = \mathop {\max }\limits_k \hat\gamma^{k}$, $\varepsilon_{A}= \mathop {\max }\limits_k \hat\varepsilon_{A}^{k}$, and PBT represents the physical boundary terms \cite{Kopriva:2016il},
\begin{equation}
\begin{split}
\PBT&= 2\sum\limits_{Boundary\atop Faces} {\int_{\partial E,N} {{{\left( {{\tilde{\statevec F}^*} - \frac{1}{2}\statevec F \cdot \hat n} \right)}^T}\statevec UdS} }
\\&=\sum\limits_{Boundary\atop Faces} {\int_{\partial E,N} {\left\{\statevec U^{T}\tilde{\mathcal A}^{+}\statevec U
+\left(\statevec U-\statevec g\right)^{T}\left|\tilde{\mathcal A}^{-}\right|\left(\statevec U-\statevec g\right)
-\statevec g^{T}\left|\tilde{\mathcal A}^{-}\right|\statevec g\right\}dS} }.
\end{split}
\label{eq:PBTGeneral}
\end{equation}
The middle term in the PBT represents additional damping due to the weak imposition of the boundary conditions through the numerical flux.
The interior face contributions satisfy
\cite{Gassner:2013uq}, 
\begin{equation}
{\tilde{\mathbf{F}}^{*T}}\jump{ {{\mathbf{U}}} } - \frac{1}{2}\jump{ {{\tilde{\mathbf{F}}^T}{\mathbf{U}}} } = -\frac{1 }{2}{\jump{ {{\mathbf{U}}} }^T}\left| \tilde{\mathcal A}\cdot\hat n \right|\jump{ {{\mathbf{U}}} }\le 0,
\label{eq:BndryDissip}
\end{equation}
pointwise.

If we define the discrete broken norm
\begin{equation}
\left\| \statevec U \right\|_N^2 = \sum\limits_{k = 1}^K {\left\| {{\statevec U^k}} \right\|_{J,N}^2},
\end{equation}
then the time derivative of the energy over the whole domain satisfies
\begin{equation}
\frac{d}{dt}\inorm{\statevec U}^{2}_{N}\le -PBT-  \sum_{\interiorfaces}\int_{\partial E,N }{\jump{ {{\mathbf{U}}} }^T}\left| \tilde{\mathcal A}\cdot\hat n \right|\jump{ {{\mathbf{U}}} }dS+2\left(\hat \gamma+\varepsilon_{a}\right)\inormN{\statevec U}^{2}.
\label{eq:GlobalStandardDiscreteNRGTimeDeriv}
\end{equation}
Equation (\ref{eq:GlobalStandardDiscreteNRGTimeDeriv}) matches that of the continuous analysis, (\ref{eq:TrueenergyTimeDerivative}), except for the product rule error term. If no energy is added through the boundaries, i.e. if $\statevec g = 0$, and if $\hat \gamma = 0$ when the same is true of the continuous problem, then 
\begin{equation}
\begin{split}
\frac{d}{dt}\inorm{\statevec U}^{2}_{N}\le &-\sum\limits_{Boundary\atop Faces} {\int_{\partial E,N} {\left\{\statevec U^{T}\tilde{\mmatrix A}^{+}\statevec U
+\statevec U^{T}\left|\tilde{\mmatrix A}^{-}\right|\statevec U
\right\}dS} }\\&-  \sum_{\interiorfaces}\int_{\partial E,N }{\jump{ {{\mathbf{U}}} }^T}\left| \tilde{\mmatrix A} \right|\jump{ {{\mathbf{U}}} }dS\\&+2\varepsilon_{A}\inormN{\statevec U}^{2}.
\end{split}
\label{eq:GlobalStandardDiscreteNRGTimeDeriv2}
\end{equation}
We see in \eqref{eq:GlobalStandardDiscreteNRGTimeDeriv2} that if the dissipation due to the internal face jumps and the weak imposition of the physical boundary conditions are not sufficiently large, then the approximation can grow in time with a growth rate, $\varepsilon_{A}$, that depends on the size of the aliasing error.

\subsection{Fully Overintegrated: $L=M,M>N$}
\label{sec:FullyOverIntegrated}
The idea behind overintegration is to approximate the integrals by quadrature of order $M>N$ sufficiently large so as to get a stable approximation. The formal statement of the weak form fully overintegrated approximation is
\begin{equation}
 \iprodN{J{\mathbf{U}}_{t},\boldsymbol\phi } + \int_{\partial E,M} {{{\tilde {\mathbf{F}} }^{*,T}}\boldsymbol\phi dS} - {\iprod{  \contravec {\statevec F},{\nabla  }\boldsymbol\phi }_M}= 0.
 \label{eq:DGSEM-S}
\end{equation}
Taking $\boldsymbol\phi=\statevec U$,
\begin{equation}
 \frac{1}{2}\frac{d}{{dt}}\left\| U \right\|_{J,N}^2 +\int_{\partial E,M} {{{\tilde {\mathbf{F}} }^{*,T}}\statevec U dS} - {\iprod{  \contravec {\statevec F},{\nabla  }\statevec U}_M}= 0.
 \label{eq:DGSEM-SE}
\end{equation}
Using the discrete extended Gauss law gives the strong form
\begin{equation}
\frac{1}{2}\frac{d}{dt}\inorm{\statevec U}^{2}_{J,N}+\int_{\partial E, M}{\statevec U^{T}\left\{\contravec{\statevec F}^{*} -\contravec{\statevec F}\cdot\hat n\right\} dS}+\iprod{\nabla\cdot\contravec{\statevec F},\statevec U}_{M} = 0.
\label{eq:WeakContinuousFormNMMStrong}
\end{equation}

In \cite{Kirby2003} it was suggested to find $M$ sufficiently large so that the quadrature error is eliminated.
We note that if $\contravec{\mathcal A}\in\mathbb{P}^{N}$ (see \eqref{eq:CoeffMatrixPoly}), then $\left(\nabla\cdot\contravec {\statevec F}\right)^{T}\statevec U = \left(\nabla\cdot\left(\contravec{\mathcal A}\statevec U\right)\right)^{T}\statevec U\in\mathbb{P}^{3N}$ so that 
\begin{equation}
\iprod{\nabla\cdot\contravec{\statevec F},\statevec U}_{M} = \iprod{\nabla\cdot\contravec{\statevec F},\statevec U}\end{equation}
if $2M-1 = 3N$, i.e., if $M>3N/2$. Alternatively, with the flux approximation \eqref{eq:fInP3N}, $2M-1=4N$ and we must take $M>2N$ for the quadrature to be exact.

The amount by which the quadratures need to be overintegrated for \emph{stability} is greater than that  needed to eliminate just the quadrature error. As long as $M\ge N$, the inner product on the right in \eqref{eq:WeakContinuousFormNMMStrong} satisfies the summation by parts property \cite{gassner2010} and so 
\begin{equation}
\iprod{\mathbb{I}^{M} \left(\contravec{\mathcal A}\statevec U\right),\nabla \statevec U }_{M}=\int_{\partial E,M}{\statevec U^{T} \mathbb{I}^{M}\left(\contravec{\mathcal A}\cdot\hat n\statevec U\right)dS} - \iprod{\nabla\cdot\mathbb{I}^{M} \left(\contravec{\mathcal A}\statevec U\right),\statevec U}_{M}.
\end{equation}
However, the crucial step needed to prove stability is to commute the differentiation and interpolation in $\nabla\cdot\mathbb{I}^{M} \left(\contravec{\mathcal A}\statevec U\right)$ 
and use the product rule \cite{Kopriva:2017yg}. The product rule holds if the interpolation is exact, i.e. 
$\mathbb{I}^{M} \left(\contravec{\mathcal A}\statevec U\right) = \contravec{\mathcal A}\statevec U$, 
which occurs if $M\ge 2N$ if $\tilde A\in\mathbb{P}^{N}$ and $M\ge 3N$ if $\tilde A\in\mathbb{P}^{2N}$. 

With $M$ sufficiently large so that $\contravec{\statevec F}=\mathbb{I}^{M} \left(\contravec{\mathcal A}\statevec U\right) = \contravec{\mathcal A}\statevec U$, 
\begin{equation}
 \iprod{\nabla\cdot\mathbb{I}^{M} \left(\contravec{\mathcal A}\statevec U\right),\statevec U}_{M} =  \iprod{\nabla\cdot \left(\contravec{\mathcal A}\statevec U\right),\statevec U}_{M}= \iprod{\left(\nabla\cdot\contravec{\mathcal A}\right)\statevec U+\contravec{\mathcal A}\cdot\nabla\statevec U ,\statevec U }_{M}.
\end{equation}
Therefore, (c.f. \eqref{eq:DGLOnFlux})
\begin{equation}
\begin{split}
\iprod{\nabla\cdot\contravec{\statevec F},\statevec U}_{M} &= \frac{1}{2}\iprod{\nabla\cdot\contravec{\statevec F},\statevec U}_{M}+ \frac{1}{2}\iprod{\left(\nabla\cdot\contravec{\mathcal A}\right)\statevec U+\contravec{\mathcal A}\cdot\nabla\statevec U ,\statevec U }_{M} \\&=
\frac{1}{2}\int_{\partial E, M}{\statevec U^{T}\contravec{\statevec F}dS} + \frac{1}{2}\iprod{\left(\nabla\cdot \contravec{\mathcal A}\right)\statevec U,\statevec U}_{M}
\end{split}
\end{equation}
%
%\begin{equation}
% \frac{1}{2}\iprod{ \statevec U,\left(\nabla\cdot\contravec{\mathcal A}\right)\statevec U}_{M} = \frac{1}{2}\int_{\partial E,M} {{{\left( {\contravec{\mathcal A}\cdot \hat n\statevec U } \right)}^T}\statevec U dS}  -\iprod{\contravec{\mathcal A}\statevec U,\nabla \statevec U}_{M},
%\label{eq:OverintegratedVolumeEnergy}
%\end{equation}
and the time derivative of the local energy is
\begin{equation}
 \frac{1}{2}\frac{d}{{dt}}\left\| U \right\|_{J,N}^2 + \int_{\partial E,M} {{{\left\{ {{{\tilde {\mathbf{F}} }^*} -\frac{1}{2} {\tilde {\mathbf{F}} } \cdot \hat n} \right\}}^T}\statevec U dS}  = -\frac{1}{2}{\iprod{\statevec U,\left(\nabla  \cdot  {\contravec{\mathcal A} } \right)\statevec U }_M}.
 \label{eq:DGSEM-SE2}
\end{equation}
The term on the right of \eqref{eq:DGSEM-SE2} can be bounded using the equivalence of the discrete and continuous norms \cite{CHQZ:2006} by
 \begin{equation}
 \begin{split}
 - {\iprod{\statevec U,\left(\nabla  \cdot  {\contravec{\mathcal A}} \right)\statevec U}_M} &=  - {\iprod{\statevec U,\left(\nabla  \cdot  {\contravec{\mathcal A}} \right)\statevec U}} \leqslant C\mathop {\max }\limits_E {\left\| {\nabla  \cdot {\contravec{\mathcal A}} } \right\|_2}\iprodN{\statevec U,\statevec U} \\&\leqslant C\frac{{\mathop {\max }\limits_E {\left\| {\nabla  \cdot {\contravec{\mathcal A}} } \right\|_2}}}{{\mathop {\min }\limits_E J}}{\left( {J\statevec U,\statevec U} \right)_N} \equiv 2{{\hat \gamma }_{OI}}\left\| \statevec U \right\|_{J,N}^2.
 \end{split}
  \end{equation}
Then for element $e^{k}$,
\begin{equation}
\frac{1}{2}\frac{d}{{dt}}\left\| {\mathbf{U}^{k}} \right\|_{J,N}^2 \leqslant  - \int_{\partial E,M} {{{\left\{ {{{\contravec {\mathbf{F}} }^*} - \frac{1}{2}\tilde {\mathbf F}^{k} \cdot \hat n} \right\}}^T}{\mathbf{U}}^{k}dS}  + \frac{1}{2}2\hat \gamma_{OI}^{k} \left\| {\mathbf{U}}^{k} \right\|_{J,N}^2.
\label{eq:OITimeDerivativeBound1}
\end{equation}
The surface fluxes are of the same form as for the Standard Approximation, see \eqref{eq:LocalDiscreteEnergyDGSEM}. Therefore, when we sum over all of the elements
\begin{equation}
\frac{d}{dt}\inorm{\statevec U}^{2}_{N}\le -PBT- BI+2\hat\gamma_{OI}\inormN{\statevec U}^{2},
\label{eq:GlobalOIFullDiscreteNRGTimeDeriv}
\end{equation}
where, now
\begin{equation}
\begin{split}
\PBT=\sum\limits_{Boundary\atop Faces} {\int_{\partial E,M} {\left\{\statevec U^{T}\tilde{\mathcal A}^{+}\statevec U
+\left(\statevec U-\statevec g\right)^{T}\left|\tilde{\mathcal A}^{-}\right|\left(\statevec U-\statevec g\right)
-\statevec g^{T}\left|\tilde{\mathcal A}^{-}\right|\statevec g\right\}dS} },
\end{split}
\end{equation}
$\tilde{\mathcal A}^{\pm}=\tilde{\mathcal A}\cdot\hat n \pm \left|\tilde{\mathcal A}\cdot\hat n\right|$ and
\begin{equation}
BI= \int_{\partial E,M }{\jump{ {{\mathbf{U}}} }^T}\left| \tilde{\mathcal A} \right|\jump{ {{\mathbf{U}}} }dS\ge 0,
\end{equation}
again, provided that the interpolant of the coefficient matrix is continuous across element interfaces.

In the special case where the energy should not grow, i.e. where $\statevec g=0$ and $\nabla  \cdot \tilde{\mathcal A}=0$, the evolution of the energy of the solution is governed by 
\begin{equation}
\frac{d}{dt}\inorm{\statevec U}^{2}_{N}\le -\sum\limits_{Boundary\atop Faces} {\int_{\partial E,M} {\left\{\statevec U^{T}\tilde{\mathcal A}^{+}\statevec U
+\statevec U^{T}\left|\tilde{\mathcal A}^{-}\right|\statevec U
\right\}dS} }-  \sum_{\interiorfaces}\int_{\partial E,M }{\jump{ {{\mathbf{U}}} }^T}\left| \tilde{\mathcal A} \right|\jump{ {{\mathbf{U}}} }dS\le 0
\label{eq:GlobalOIFullDiscreteNRGTimeDeriv2}
\end{equation}
and does not grow in time, implying stability.

In summary, when the flux is a polynomial of its argument, and the volume and the surface integrals are approximated to a sufficiently high order, $L = M > 2N$ for $\contravec{\mathcal A}\inPN$ and  $L = M > 3N$ for $\contravec{\mathcal A}\in\mathbb{P}^{2N}$, then the product rule error seen in standard, underintegrated approximation, does not appear, and the overintegrated approximation is stable. 

\subsection{Volume Only OverIntegrated: $L=N, M>N$}
As was originally suggested, overintegration could be performed to eliminate the product rule error term $\varepsilon_{A}$ in the volume, but not along the faces \cite{Mengaldo201556}. In the partially overintegrated approximation, we can construct two more ``strong'' forms, one of which differs from the ``weak'' form. \emph{In the following, we will assume that $\statevec g=0$ and ${\nabla  \cdot{ \contravec{\mathcal A} }}=0$ so that any growth in the energy represents instability.}

The strong and weak forms of the approximations are no longer necessarily algebraically equivalent, so we must analyze the two separately.
The weak form approximation with the surface terms underintegrated is
\begin{equation}
[W]\quad\iprod{ J\statevec U,\phi}_{N}+\int_{\partial E, N}{\phi^{T}\contravec{\statevec F}^{*}dS}-\iprod{\contravec{\statevec F},\nabla\phi}_{M} = 0.
\label{eq:WeakContinuousFormNNM}
\end{equation}
We can also construct two ``strong'' form approximations. We can integrate by parts first and then apply quadrature. Starting from (\ref{eq:approxWExactIntegration}), applying the Gauss law and then replacing integrals by quadrature gives the first strong form
\begin{equation}
[S1]\quad\left\langle {{J{\mathbf{U}}_t},\phi } \right\rangle  + \int_{\partial E,N} {{\mathbb{I}^N}\left( {{{\left( {{\tilde{\mathbf{F}}^*} - \tilde {\mathbf{F}}  \cdot \hat n} \right)}^T}\phi } \right)dS}  + {\left\langle {\nabla  \cdot {\mathbf{\spacevec F}},\phi } \right\rangle _M} = 0.
\label{eq:OIStrong1}
\end{equation}
As an aside, form [S1] is not particularly interesting in practice as it is not conservative if $N\ne M$.
The second form we get by applying the discrete Gauss Law to the weak form, (\ref{eq:WeakContinuousFormNNM}), giving
\begin{equation}
[S2]\quad{\left\langle {{J{\mathbf{U}}_t},\phi } \right\rangle _N} + \int_{\partial E,N} {{\mathbb{I}^N}\left( {{\tilde{\mathbf{F}}^{*,T}}\phi } \right)dS}  - \int_{\partial E,M} {{{\mathbf{F}}^{*,T}}\phi dS}  + {\left\langle {\nabla  \cdot {\mathbf{\spacevec F}},\phi } \right\rangle _M} = 0
\label{eq:OIStrongForm2}
\end{equation}
The approximations (\ref{eq:OIStrong1}) and (\ref{eq:OIStrongForm2}) are not the same unless $N=M$. Note that (\ref{eq:OIStrongForm2}) is algebraically equivalent to the weak form (\ref{eq:WeakContinuousFormNNM}), $[W]\Leftrightarrow[S2]$,  (c.f. \cite{gassner2010}). Therefore we need to analyze the stability of [S1] and one of $[W]$ or $[S2]$.

\subsubsection{Stability of the Approximation [S1]}
To find the stability property of $[S1]$, we replace $\phi$ by $\statevec U$, as usual. The volume term can then be replaced as in Sec. \ref{sec:FullyOverIntegrated} so that the elemental energy is governed by
\begin{equation}\frac{1}{2}\frac{d}{{dt}}{\left\| {\mathbf{U}} \right\|^2_{J,N}} + \int_{\partial E,N} {{\mathbb{I}^N}\left( {{{\left( {{{\mathbf{F}}^*} - \tilde {\mathbf{F}}  \cdot \hat n} \right)}^T}\mathbf U } \right)dS}  + \frac{1}{2}\int_{\partial E,M} {{{\left( {\tilde {\mathbf{F}}  \cdot \hat n} \right)}^T}{\mathbf{U}}dS}  = 0.
\label{eq:S1Energy1}
\end{equation}
We can add and subtract terms to see the contribution from the inconsistent integration. Since
$M$ is chosen so that the quadrature is exact, \eqref{eq:S1Energy1} is equivalent to
\begin{equation}
\begin{split}
\frac{1}{2}\frac{d}{{dt}}{\left\| {\mathbf{U}} \right\|^2_{J,N}} &+ \int_{\partial E,N} {{\mathbb{I}^N}\left( {{{\left( {{{\mathbf{F}}^*} - \frac{1}{2}\tilde {\mathbf{F}}  \cdot \hat n} \right)}^T}\mathbf U } \right)dS}  \\&+ \frac{1}{2}\int_{\partial E} {\left\{{{\left( {\tilde {\mathbf{F}}  \cdot \hat n} \right)}^T}{\mathbf{U}} - \mathbb{I}^{N}\left( {{\left( {\tilde {\mathbf{F}}  \cdot \hat n} \right)}^T}{\mathbf{U}} \right)\right\}dS}  = 0.
\end{split}
\label{eq:S1Energy2}
\end{equation}

Summing over all of the elements gives
\begin{equation}\frac{d}{{dt}}\sum\limits_{k = 1}^K {{{\left\| {{{\mathbf{U}}^k}} \right\|}^2_{J,N}}} = -BI - PBT,
\label{eq:UnderIntEnergy1}
\end{equation}
where the  interior interface contributions are now
\begin{equation}
\begin{split}
BI = &-2\sum\limits_{\interiorfaces} \int_{\partial E,N} {\left\{ {{{\mathbf{F}}^{*,T}}\jump{\mathbf{U}} - \frac{1}{2}\jump{{{\left( {\tilde {\mathbf{F}}  \cdot \hat n} \right)}^T}{\mathbf{U}}}} \right\}dS}  \\&- \sum\limits_{\interiorfaces}\int_{\partial E} {\left\{\jump{{{\left( {\tilde {\mathbf{F}}  \cdot \hat n} \right)}^T}{\mathbf{U}}} - \jump{\mathbb{I}^{N}\left({{\left( {\tilde {\mathbf{F}}  \cdot \hat n} \right)}^T}{\mathbf{U}}\right)}\right\}dS} 
\end{split}
\label{eq:InterfaceBoundaryTerms}
 \end{equation}
and the physical boundary terms are now
\begin{equation}
\begin{split}
PBT = 2&\sum\limits_{\boundaryfaces} \int_{\partial E,N} { {{{\left( {{{\mathbf{F}}^*} - \frac{1}{2}{{\tilde {\mathbf{F}} }} \cdot \hat n} \right)}^T}{{\mathbf{U}}}} dS}  \\&+\sum\limits_{\boundaryfaces}\int_{\partial E} {\left\{{{\left( {{{\tilde {\mathbf{F}} }} \cdot \hat n} \right)}^T}{\mathbf{U}} -\mathbb{I}^{N}\left( {{\left( {{{\tilde {\mathbf{F}} }} \cdot \hat n} \right)}^T}{\mathbf{U}}\right)\right\}dS}  .
\end{split}
\end{equation}
 
 We examine the interior boundary contributions first.  From the algebraic identity,
\begin{equation}\jump{ {\left( {\tilde {\mathbf{F}}  \cdot \hat n} \right)^T}{\mathbf{U}}}  = \average{ {\left( {\tilde {\mathbf{F}}  \cdot \hat n} \right)^T}} \jump{ {\mathbf{U}}}  + \jump{ {\left( {\tilde {\mathbf{F}}  \cdot \hat n} \right)^T}} \average{ {\mathbf{U}}}, \end{equation}
and the particular linear structure of the flux,
\begin{equation}\jump{ {\left( {\tilde {\mathbf{F}}  \cdot \hat n} \right)}{\mathbf{U}}}^T  = 2\average{{ {\mathbf{U}}}} ^T\tilde{\mathcal A} \jump{ {\mathbf{U}}} \end{equation}
at each surface point. Then using \eqref{eq:BndryDissip}, and the fact that the interpolant of the jumps is equal to the jump of the interpolants, we see that the interior interfaces have a dissipative component and an interpolation error component whose sign is indeterminate
\begin{equation}
\begin{split}
BI = &\sum\limits_{\interiorfaces} \int_{\partial E,N} {\jump{ {{\mathbf{U}}} }^T}\left| \tilde{\mathcal A}\cdot\hat n \right|\jump{ {{\mathbf{U}}} }dS  \\&- \sum\limits_{\interiorfaces}\int_{\partial E} {\left\{\average{{ {\mathbf{U}}}} ^T\tilde{\mathcal A}\cdot\hat n \jump{ {\mathbf{U}}} -\mathbb{I}^{N}\left(\average{{ {\mathbf{U}}}} ^T\tilde{\mathcal A}\cdot\hat n \jump{ {\mathbf{U}}}\right) \right\}dS} .
\end{split}
\label{eq:InterfaceBoundaryTerms2}
 \end{equation}

A similar error is introduced at the physical boundaries. 
%Let us add and subtract terms to rewrite $PBT$ as
%\begin{equation}
%\begin{split}
%PBT =  &- \sum\limits_{{{ \boundaryfaces}}} \int_{\partial E,N} {{\mathbb{I}^N}\left( {{{\left( {{{\mathbf{F}}^*} - \frac{1}{2}{{\widetilde {\mathbf{F}}}^k} \cdot \hat n} \right)}^T}{{\mathbf{U}}^k}} \right)dS}  \\&+ \frac{1}{2}\sum\limits_{{{ \boundaryfaces}}}\left\{ {\int_{\partial E} {{{\left( {{{\widetilde {\mathbf{F}}}^k} \cdot \hat n} \right)}^T}{{\mathbf{U}}^k}dS - \int_{\partial E,N} {{\mathbb{I}^N}\left( {{{\left( {{{\widetilde {\mathbf{F}}}^k} \cdot \hat n} \right)}^T}{{\mathbf{U}}^k}} \right)dS} } } \right\}
%\end{split}
%\label{eq:PBTUnderInt1}
% \end{equation}
 Following \eqref{eq:PBTGeneral},
 \begin{equation}
\begin{split}
PBT = &\sum\limits_{\boundaryfaces} \int_{\partial E,N} { \left({{{\mathbf{U}}^T}{{\tilde {\mathcal{A}}}^ + }{\mathbf{U}} + {{\mathbf{U}}^T}\left| {{{\tilde{\mathcal{A}}}^ - }} \right|{\mathbf{U}}}\right) dS} 
 \\&+\sum\limits_{\boundaryfaces}\int_{\partial E} {\left\{{\statevec U^T}\tilde{\mathcal A}\cdot\hat n{\mathbf{U}} -\mathbb{I}^{N}\left({\statevec U^T}{\tilde {\mmatrix{A}}}\cdot\hat n{\mathbf{U}} \right)\right\}dS}  .
\end{split}
\label{eq:PBT_Underintegrated}
\end{equation}
When we insert \eqref{eq:InterfaceBoundaryTerms2} and \eqref{eq:PBT_Underintegrated} into \eqref{eq:UnderIntEnergy1},
\begin{equation}
\begin{split}
\frac{d}{dt}\inormN{\statevec U} ^{2}= &-\sum\limits_{\interiorfaces} \int_{\partial E,N} {\jump{ {{\mathbf{U}}} }^T}\left| \tilde{\mathcal A} \cdot\hat n\right|\jump{ {{\mathbf{U}}} }dS 
 \\&+ \sum\limits_{\interiorfaces}\int_{\partial E} {\left\{\average{{ {\mathbf{U}}}} ^T\tilde{\mathcal A}\cdot\hat n \jump{ {\mathbf{U}}} -\mathbb{I}^{N}\left(\average{{ {\mathbf{U}}}} ^T\tilde{\mathcal A}\cdot\hat n \jump{ {\mathbf{U}}}\right) \right\}dS}\\&- \sum\limits_{\boundaryfaces} \int_{\partial E,N} { \left({{{\mathbf{U}}^T}{{\tilde {\mathcal{A}}}^ + }{\mathbf{U}} + {{\mathbf{U}}^T}\left| {{{\tilde{\mathcal{A}}}^ - }} \right|{\mathbf{U}}}\right) dS} 
\\& -\sum\limits_{\boundaryfaces}\int_{\partial E} {\left\{{\statevec U^T}\tilde{\mathcal A}\cdot\hat n{\mathbf{U}} -\mathbb{I}^{N}\left({\statevec U^T}{\tilde {\mathcal{A}}}\cdot\hat n{\mathbf{U}} \right)\right\}dS}
 \end{split}
 \label{eq:NNMTotalEnergyEqn}
\end{equation}
We see that the approximation has element face dissipation plus interpolation error terms of indeterminate sign that come from the fact that the inconsistent quadrature does not represent the surface fluxes exactly. If that interpolation error is large, i.e. in severely under resolved approximations,  then it is possible that the right hand side of \eqref{eq:NNMTotalEnergyEqn} is positive and the energy of the solution grows in time.

% \begin{equation}{\mathbb{I}^N}\left( {{{\mathbf{U}}^T}{{\tilde {\mathcal{A}}}^ + }{\mathbf{U}} + {{\mathbf{U}}^T}\left| {{{\tilde{\mathcal{A}}}^ - }} \right|{\mathbf{U}}} \right) \geqslant 0.\end{equation}
% Then on each face,
% \begin{equation}PB{T^i} =  - \int_{\partial {E^i},N} {{\mathbb{I}^N}\left( {{{\mathbf{U}}^T}{{\tilde{ \mathcal{A}}}^ + }{\mathbf{U}} + {{\mathbf{U}}^T}\left| {{{\tilde {\mathcal{A}}}^ - }} \right|{\mathbf{U}}} \right)dS}  + \frac{1}{2}\int_{\partial {E^i}} {\left\{ {{{\left( {\tilde {\mathbf{F}} \cdot \hat n} \right)}^T}{\mathbf{U}} - {\mathbb{I}^N}\left( {{{\left( {\tilde {\mathbf{F}} \cdot \hat n} \right)}^T}{\mathbf{U}}} \right)} \right\}dS} \end{equation}
%As in the interior, the physical boundary terms have dissipation plus a sign indeterminate term that depends on the size of the interpolation error, which could lead to energy growth. Therefore,
%\begin{equation}
%\frac{d}{dt}\inorm{\statevec U}^{2}_{J,N}\le\sum_{{\boundaryfaces}}
%\end{equation}

\subsubsection{Stability of the Approximation $[S2]\Leftrightarrow [W]$}

Again, we replace $\phi$ with $\statevec U$, this time in \eqref{eq:WeakContinuousFormNNM}. 
Then
\begin{equation}\frac{1}{2}\frac{d}{{dt}}\left\| \statevec U \right\|_{J,N}^2 + \int_{\partial E,N} {{\statevec U^T}{{\tilde {\statevec F}}^*}dS}  - \frac{1}{2}\int_{\partial E,M} {{\statevec U^T}\tilde {\statevec F} \cdot \hat ndS}  = 0.
\label{eq:S2LocalEnergy}\end{equation}
We add and subtract terms and use exactness of the consistently integrated quadrature to rewrite \eqref{eq:S2LocalEnergy} as
\begin{equation}\frac{1}{2}\frac{d}{{dt}}\left\| \statevec U \right\|_{J,N}^2 + \int_{\partial E,N} {{\statevec U^T}\left( {{{\tilde {\statevec F}}^*} - \frac{1}{2}\tilde {\statevec F} \cdot \hat n} \right)dS}  - \frac{1}{2}\int_{\partial E} {\left\{ {{\statevec U^T}\tilde {\statevec F} \cdot \hat n - {\mathbb{I}^N}\left( {{\statevec U^T}\tilde {\statevec F} \cdot \hat n} \right)} \right\}dS}  = 0.
\label{eq:S2Energy}
\end{equation}
The only difference between \eqref{eq:S1Energy2} and \eqref{eq:S2Energy} is in the sign of the interpolation error term. Therefore the
energy analysis of the previous section applies to the second strong form approximation in the same way.

\subsubsection{Analysis of the Element Boundary Contributions}

We first show that the errors due to inconsistent integration of the surface terms are due to aliasing of the modes $2N,\ldots,pN$, where $p=3$ or $p=4$, depending on how the coefficient matrices are approximated.  The interpolation errors in \eqref{eq:NNMTotalEnergyEqn} are of the form
\begin{equation}
\int_{\partial E}{\left\{V - \mathbb{I}^{N}(V)\right\}dS},
\end{equation}
where $V\in\mathbb{P}^{pN}$. To simplify the analysis, we will consider two dimensional geometries where the element faces
are the element boundary curves. Then along each side, the surface integral reduces to
\begin{equation}
\epsilon \equiv \int_{-1}^{1}\left\{V(\xi)- \mathbb{I}^{N}(V(\xi))\right\}d\xi=\iprod{V-\mathbb{I}^{N}(V),1}=\iprod{V-\mathbb{I}^{N}(V),L_{0}},
\label{eq:InterpErrorIntegral}
\end{equation}
where $L_{0}=1$ is the Legendre polynomial of degree zero.
In modal form, 
\begin{equation}
V = \sum_{k=0}^{pN}{\hat V_{k} L_{k}(\xi)}
\end{equation}
where $\hat V_{k}=\iprod{V,L_{k}}/\inorm{L_{k}}^{2},\;k=0,1,\ldots,pN$ are the modal coefficients of the polynomial and $L_{k}$ is the Legendre polynomial of degree $k$. The modal representation of the interpolant of $V$ is
\begin{equation}
\mathbb{I}^{N}(V) = \sum_{k=0}^{N}{\bar V_{k} L_{k}(\xi)},
\end{equation}
where \cite{CHQZ:2006}
\begin{equation}{\bar V_k} = {\hat V_k} + a_{k},\end{equation}
and the $a_{k}$ are the aliases of the true coefficients
\begin{equation}
a_{k}=\frac{1}{{\left\| {{L_k}} \right\|_N^2}}\sum\limits_{n = N + 1}^{pN} {{{\left\langle {{L_n},{L_k}} \right\rangle }_N}{{\hat V}_n}}.
\end{equation}
Then
\begin{equation}V - {\mathbb{I}^N}(V) = \sum\limits_{k = N + 1}^{pN} {{{\hat V}_k}{L_k}}  - \sum\limits_{k = 0}^N {{a_k}{L_k}}. \end{equation}
Orthogonality of the Legendre polynomials causes most of the terms in \eqref{eq:InterpErrorIntegral} to vanish, leaving
\begin{equation} \epsilon=  - 2{a_0} =  - \sum\limits_{n = N + 1}^{pN} {{{\left\langle {{L_n},{L_0}} \right\rangle }_N}{{\hat V}_n}} .\end{equation}
Finally, $\iprodN{L_{n},L_{0}} = \iprod{L_{n},L_{0}}=0$ for $n=N+1,N+2,\ldots,2N-1$ so
\begin{equation}
\int_{-1}^{1}\left\{V(\xi)- \mathbb{I}^{N}(V(\xi))\right\}d\xi= - \sum\limits_{n = 2N}^{pN} {{{\left\langle {{L_n},{L_0}} \right\rangle }_N}{{\hat V}_n}}.
\end{equation}
There is no error if, in fact, $V\in\mathbb{P}^{2N-1}$. Otherwise, the error is due to aliasing of the modes from $2N$ to $pN$.

It remains, then, to see if the aliasing error is smaller or larger than the dissipation due to the numerical flux. First we note that
the internal interface dissipation in \eqref{eq:NNMTotalEnergyEqn} is $O\left( \jump{U}^{2}\right)$, whereas the quadrature aliasing error
is $O\left(\jump{U} \right)$. In lieu of finding the coefficients of a triple product, we illustrate the dissipation and the aliasing errors along an interior edge
\begin{equation}
-\int_{-1,N}^{1}{\jump{ {{\mathbf{U}}} }^T}\left| \tilde{\mathcal A} \cdot\hat n\right|\jump{ {{\mathbf{U}}} }d\xi + \int_{-1}^{1} {\left\{\average{{ {\mathbf{U}}}} ^T\tilde{\mathcal A} \cdot\hat n\jump{ {\mathbf{U}}} -\mathbb{I}^{N}\left(\average{{ {\mathbf{U}}}} ^T\tilde{\mathcal A} \cdot\hat n\jump{ {\mathbf{U}}}\right) \right\}d\xi}
\end{equation}
for a scalar problem with
\begin{equation}\begin{gathered}
  \jump{ U}  = \alpha {\left( {1 + \xi } \right)^{q/3}} \hfill \\
  \average{ U}  = \beta {\left( {1 + \xi } \right)^{q/3}} \hfill \\
  {\tilde {\mathcal{A}}}\cdot\hat n = \gamma {\left( {1 + \xi } \right)^{q/3}}. \hfill \\ 
\end{gathered} \end{equation} 
With this choice, an interior edge contribution to the energy is
\begin{equation} - {\alpha ^2}\left| \gamma  \right|\int_{ - 1,N}^1 {{{\left( {1 + \xi } \right)}^q}d\xi }  + \alpha \beta \gamma \int_{ - 1}^1 {\left\{ {{{\left( {1 + \xi } \right)}^q} - {\mathbb{I}^N}\left( {{{\left( {1 + \xi } \right)}^q}} \right)} \right\}}. \end{equation}
As noted above, if $q\le 2N-1$ the aliasing contribution should vanish and the edge contribution is dissipative. 

If $q$ is large, then the edge contributions can be destabilizing when the aliasing errors are
large enough, i.e. in underresolved problems. As a concrete example, we choose $\alpha = 10^{-3}$, $\beta = 1$ and $\gamma = -1$. Table \ref{tab:ErrorTable} shows the edge dissipation, aliasing contribution and total dissipation for $q = 18$ and approximation orders $N$ between three and 16. We see that for $N\le 5$, the aliasing error dominates the dissipation due to upwinding and the combined contribution is positive. As the resolution increases, the aliasing error decays exponentially fast (as seen in a semi-log plot of that term vs $N$) so that  the overall contribution is dissipative for $6\le N\le 17$. Finally, for $2N-1\ge q$ the aliasing error vanishes as expected and the only edge contribution to the energy comes from the upwinding. Table \ref{tab:ErrorTable} is therefore an illustration that severe underresolution can lead to destabilizing aliasing errors that might not be significant as the resolution increases.

\begin{table}[htp]
\begin{center}\begin{tabular}{cccccc}N & 2N-1 &$ - {\alpha ^2}\left| \gamma  \right|\int_{ - 1,N}^1 {{{\left( {1 + \xi } \right)}^q}d\xi }$ & $ \alpha \beta \gamma \int_{ - 1}^1 {\left\{ {{{\left( {1 + \xi } \right)}^q} - {\mathbb{I}^N}\left( {{{\left( {1 + \xi } \right)}^q}} \right)} \right\}}d\xi$ & Sum \\
 3  & 5   &-4.434E-02   & 1.674E+01   & {\bf 1.670E+01}\\
 4 & 7 & -3.092E-02    &3.327E+00    &{\bf 3.296E+00}\\
 5  &9 &  -2.799E-02   & 3.967E-01    &{\bf 3.687E-01}\\
 6 &11 &  -2.762E-02  &  2.560E-02   &-2.018E-03\\
 7& 13 &  -2.759E-02  &  7.737E-04   &-2.682E-02\\
 8 &15 &  -2.759E-02  &  8.237E-06  & -2.759E-02\\
 9 &17 &  -2.759E-02  &  1.297E-08   &-2.759E-02\\ \hline
10 &19 &  -2.759E-02  & -2.310E-15  & -2.759E-02\\
11 &21 &  -2.759E-02  &  1.102E-14   &-2.759E-02\\
12 &23 &  -2.759E-02  &  1.826E-14  & -2.759E-02\\
13 &25 &  -2.759E-02  & -2.218E-14  & -2.759E-02\\
\end{tabular} \caption{Interior interface dissipation along a single edge as a function of polynomial order, $N$, for $q=18$. Bold entries predict instability. The horizontal line marks where the aliasing error associated with the underintegration of the boundary terms vanishes.}
\end{center}
\label{tab:ErrorTable}
\end{table}

\section{Conclusions}

To summarize the results, let us gather the boundary and interface dissipation terms as
\[
Dissip\equiv \sum_{\interiorfaces}\int_{\partial E,N }{\jump{ {{\mathbf{U}}} }^T}\left| \tilde{\mathcal A}\cdot\hat n \right|\jump{ {{\mathbf{U}}} }dS+\sum\limits_{Boundary\atop Faces} {\int_{\partial E,N} {\left\{\statevec U^{T}\tilde{\mathcal A}^{+}\statevec U
+\statevec U^{T}\left|\tilde{\mathcal A}^{-}\right|\statevec U
\right\}dS} }\ge 0.
\]
Then for linear problems where there should be no energy growth, the energy of the standard, the partially (volume only) overintegrated, and the fully overintegrated approximations satisfy
\begin{equation}
\text{[Standard]}\quad\frac{d}{dt}\inorm{\statevec U}^{2}_{N}\le -Dissip+2\varepsilon_{a}\inormN{\statevec U}^{2},
\label{eq:GlobalStandardDiscreteNRGTimeDeriv2S}
\end{equation}
\begin{equation}
\begin{split}
\text{[Volume Overintegrated]} \quad\frac{d}{dt}\inormN{\statevec U} ^{2}= 
&-Dissip 
 \\&\pm \sum\limits_{\interiorfaces}\int_{\partial E} {\left\{\average{{ {\mathbf{U}}}} ^T\tilde{\mathcal A} \cdot\hat n\jump{ {\mathbf{U}}} -\mathbb{I}^{N}\left(\average{{ {\mathbf{U}}}} ^T\tilde{\mathcal A} \cdot\hat n\jump{ {\mathbf{U}}}\right) \right\}dS}
 \\&\pm\sum\limits_{\boundaryfaces}\int_{\partial E} {\left\{{\statevec U^T}\tilde{\mathcal A}\cdot\hat n{\mathbf{U}} -\mathbb{I}^{N}\left({\statevec U^T}{\tilde {\mathcal{A}}}\cdot\hat n{\mathbf{U}} \right)\right\}dS},
 \end{split}
 \label{eq:NNMTotalEnergyEqnS}
\end{equation}
where the $\pm$ depends on which strong form is used, and
\begin{equation}
\begin{split}
\text{[Fully Overintegrated]} \quad\frac{d}{dt}\inorm{\statevec U}^{2}_{N}\le -Dissip\le 0
\end{split}
\label{eq:GlobalOIFullDiscreteNRGTimeDeriv2S}
\end{equation}
provided that $M$ is sufficiently large so that the product rule holds.

The results support the finding of \cite{Mengaldo201556} that overintegrating the surface and volume integrals leads to a more robust approximation.
The standard approximation has a growth term whose growth rate factor, $\varepsilon_{a}$, depends on aliasing errors associated with the amount by which the product rule fails to hold for the interpolants of products of polynomials. That term can lead to exponential growth if the dissipation terms associated with the boundary and interface conditions are not sufficiently large, which, in a sense, is a definition of ``severely underresolved''. The volume only overintegrated approximation eliminates that aliasing term associated with the volume, so there is no $\varepsilon_{a}$ term, but it introduces aliasing terms along the element faces. In severely underresolved problems, these aliasing terms could also be large enough to destabilize the approximation. Using consistent integration for both the surfaces and the volume leads to a stable approximation where the energy of the solution is nonincreasing. 

As a final comment, we note that none of these approximations are dealiased, yet some are stable. There is aliasing in the discrete norm of the solution stemming from the fact that the argument of $\inormN{J\statevec U}$ will be a polynomial of degree $3N$, which cannot be approximated exactly with the Gauss quadratures of order $N$. Nevertheless, the discrete norm is equivalent to the continuous norm so discrete stability implies stability in the continuous norm. Also, unless the coefficient matrices $\mmatrix A_{m}$ are polynomial functions of their argument, there will always be aliasing errors when representing the flux as a polynomial of any finite degree. Nevertheless, the fully overintegrated approximation is stable (nonincreasing energy when expected) provided that the interpolant of the coefficient matrices is constructed to be divergence free. Aliasing errors are still found, nonetheless, in the dissipation boundary terms in the approximations of the coefficient matrices, but those only affect the rate of dissipation, not stability per se. Therefore, it should be emphasized that the presence of aliasing in an approximation does not imply instability, nor does stability imply no aliasing, and that the ``fully overintegrated'' approximation has aliasing errors that do not contribute to instability.

{\bf Acknowledgement.} The author would like to thank Gregor Gassner for helpful discussions. This work was supported by a grant from the Simons Foundation (\#426393, David Kopriva).
%In that case, the energy satisfies \eqref{eq:GlobalOIFullDiscreteNRGTimeDeriv}, but with
%\begin{equation}
%\begin{split}
%\PBT=\sum\limits_{Boundary\atop Faces} {\int_{\partial E,N} {\left\{\statevec U^{T}\tilde{\mathcal A}^{+}\statevec U
%+\left(\statevec U-\statevec g\right)^{T}\left|\tilde{\mathcal A}^{-}\right|\left(\statevec U-\statevec g\right)
%-\statevec g^{T}\left|\tilde{\mathcal A}^{-}\right|\statevec g\right\}dS} }
%\end{split}
%\end{equation}
%and

\bibliographystyle{plain}

\bibliography{../../BibTex/dakBib}

\end{document}